\documentclass[UTF-8,reqno]{amsart}
\usepackage{enumerate}
\usepackage{amssymb,url,color, booktabs}
\usepackage{mathrsfs}

\usepackage[nobysame]{amsrefs}
\BibSpec{article}{%
+{}{\PrintAuthors} {author}
+{,}{ \textrm} {title}
+{.}{ \textit} {journal}
+{,}{ \textbf} {volume}
+{}{ \parenthesize} {date}
+{,}{ } {pages}
+{.}{ arXiv:} {eprint}
+{.}{} {transition}
}
\BibSpec{book}{%
+{}{\PrintAuthors} {author}
+{,}{ \textit} {title}
+{.}{ \textrm} {series} 
+{,}{ Vol.} {volume} 
+{.}{ } {publisher}
+{,}{ } {date}
+{.}{} {transition}
}
\usepackage{color}
\usepackage[colorlinks=true]{hyperref}
\hypersetup{
    linkcolor=blue,          
    citecolor=red,        
    filecolor=blue,      
    urlcolor=cyan  
}

\newcommand{\red}{\color{red}}

\numberwithin{equation}{section}

\newcommand{\be}{\begin{eqnarray}}
\newcommand{\ee}{\end{eqnarray}}
\newcommand{\ce}{\begin{eqnarray*}}
\newcommand{\de}{\end{eqnarray*}}
\newtheorem{theorem}{Theorem}[section]
\newtheorem{lemma}[theorem]{Lemma}
\newtheorem{remark}[theorem]{Remark}
\newtheorem{definition}[theorem]{Definition}
\newtheorem{proposition}[theorem]{Proposition}
\newtheorem{Examples}[theorem]{Example}
\newtheorem{corollary}[theorem]{Corollary}

\def\eps{\varepsilon}

\def\e{\mathrm{e}}

\def\p{\partial}

\def\[{{\Big[}}
\def\]{{\Big]}}
\def\<{{\langle}}
\def\>{{\rangle}}
\def\({{\Big(}}
\def\){{\Big)}}

\def\bx{{\mathbf{x}}}

\def\dif{{\mathord{{\rm d}}}}

\def\no{\nonumber}
\def\={&\!\!=\!\!&}

\def\cL{{\mathcal L}}

\def\mE{{\mathbb E}}

\def\mN{{\mathbb N}}

\def\mP{{\mathbb P}}
\def\mQ{{\mathbb Q}}
\def\mR{{\mathbb R}}

\def\mW{{\mathbb W}}

\def\bP{{\mathbf P}}

\def\bE{{\mathbf E}}
\def\1{{\mathbf{1}}}

\def\sF{{\mathscr F}}

\def\sI{{\mathscr I}}

\def\geq{\geqslant}
\def\leq{\leqslant}
\def\ge{\geqslant}
\def\le{\leqslant}

\def\div{\mathord{{\rm div}}}

\def\eps{\varepsilon}

\def\e{\mathrm{e}}

\def\p{\partial}

\def\[{{\Big[}}
\def\]{{\Big]}}
\def\<{{\langle}}
\def\>{{\rangle}}
\def\({{\Big(}}
\def\){{\Big)}}

\def\bx{{\mathbf{x}}}

\def\dif{{\mathord{{\rm d}}}}

\def\no{\nonumber}
\def\={&\!\!=\!\!&}
\def\bt{\begin{theorem}}
\def\et{\end{theorem}}
\def\bl{\begin{lemma}}
\def\el{\end{lemma}}
\def\br{\begin{remark}}
\def\er{\end{remark}}
\def\bx{\begin{Examples}}
\def\ex{\end{Examples}}
\def\bd{\begin{definition}}
\def\ed{\end{definition}}
\def\bp{\begin{proposition}}
\def\ep{\end{proposition}}
\def\bc{\begin{corollary}}
\def\ec{\end{corollary}}

\def\geq{\geqslant}
\def\leq{\leqslant}
\def\ge{\geqslant}
\def\le{\leqslant}

\def\div{\mathord{{\rm div}}}

\def\bP{{\mathbf P}}

\def\<{\langle} \def\>{\rangle}

\def\red{\color{red}}

\def\bpf{\begin{proof}}
\def\epf{\end{proof}}

\allowdisplaybreaks

\begin{document}
	
\title[Euler scheme for density dependent SDEs]{Euler scheme for density dependent stochastic differential equations}
\date{\today}
\author{Zimo Hao, Michael R\"ockner and Xicheng Zhang}

\address{Zimo Hao:
	School of Mathematics and Statistics, Wuhan University,
Wuhan, Hubei 430072, P.R.China\\
Email: zimohao@whu.edu.cn}

\address{Michael R\"ockner:
Fakult\"at f\"ur Mathematik, Universit\"at Bielefeld,
33615, Bielefeld, Germany\\
Email: roeckner@math.uni-bielefeld.de
 }

\address{Xicheng Zhang:
School of Mathematics and Statistics, Wuhan University,
Wuhan, Hubei 430072, P.R.China\\
Email: XichengZhang@gmail.com
 }

\thanks{{\it Keywords: density dependent SDE, heat kernel, Strong solution}}

\thanks{
This work is supported by NNSFC grant of China (No. 11731009)  and the DFG through the CRC 1283 
``Taming uncertainty and profiting from randomness and low regularity in analysis, stochastics and their applications''. }

\begin{abstract}
In this paper we show the existence and uniqueness for a class of density dependent SDEs with bounded measurable drift,
where the existence part is based on Euler's approximation for density dependent SDEs
 and the uniqueness is based on the associated nonlinear Fokker-Planck equation. 
As an application, we obtain the well-posedness of a nonlinear Fokker-Planck equation.
\end{abstract}

\maketitle
\setcounter{tocdepth}{2}

\section{Introduction}

In this paper we consider the following density dependent stochastic differential equation (abbreviated as DDSDE):
\begin{align}\label{SDE}
\dif X_t=b(t,X_t,\rho_t(X_t))\dif t+\sqrt{2}\dif W_t,\ \ X_0\stackrel{(d)}{=}\nu_0,
\end{align}
where $W_t$ is a standard $d$-dimensional Brownian motion on some probability space
$(\Omega,\sF,\mP)$, $b:\mR_+\times\mR^d\times\mR\to\mR^d$ is a bounded Borel measurable vector field
and for $t>0$, 
$\rho_t(x)=\mP\circ X_t^{-1}(\dif x)/\dif x$ is the distributional density of $X_t$ with respect to the Lebesgue measure $\dif x$ on $\mR^d$,
$\nu_0$ is a probability measure over $\mR^d$.
By It\^o's formula, one sees that $\rho_t$ solves the following nonlinear Fokker-Planck equation (FPE) in the distributional sense:
 \begin{align}\label{P13}
\p_t\rho_t-\Delta\rho_t+\div(b(t,\cdot,\rho_t)\rho_t)=0,\quad \lim_{t\downarrow 0}\rho_t=\nu_0\mbox{ weakly}.
\end{align}
More precisely, for any $\varphi\in C_0^\infty(\mR^d)$,
\begin{align}\label{DPDE}
\<\rho_t, \varphi\>=\<\nu_0,\varphi\>+\int_0^t\<\rho_s,\Delta\varphi\>\dif s
+\int_0^t\<\rho_s, b(s,\cdot,\rho_s)\cdot\nabla\varphi\>\dif s,
\end{align}
where $\<\rho_t, \varphi\>:=\int_{\mR^d}\varphi(x)\rho_t(x)\dif x=\mE\varphi(X_t)$.

\medskip

Since the coefficients of SDE \eqref{SDE} depend on the distributional density of the solution $X_t$ evaluated at $X_t$, 
the dependence of $b$ on the measure $\rho_t(x)\dif x$, is called ``Nemytskii-type'' dependence (cf. \cite{Ba-Ro18a}, \cite{Ba-Ro18}). 
Thus \eqref{SDE} can be also called McKean-Vlasov SDE of Nemytskii-type. Let us first
recall the definition of a weak solution to DDSDE \eqref{SDE}:
\bd\label{Def1}
Let $\nu_0$ be a probability measure on $\mR^d$. We call a filtered probability space
$(\Omega,\sF,\bP; (\sF_t)_{t\geq 0})$ together with a pair of processes $(X, W)$ defined on it 
a  weak solution of SDE \eqref{SDE} with initial distribution $\nu_0$, if 
\begin{enumerate}[(i)]
\item $\bP\circ X_0^{-1}=\nu_0$, $W$ is a $d$-dimensional $\sF_t$-Brownian motion.
\item for each $t>0$, $\bP\circ X^{-1}_t(\dif x)/\dif x=\rho_t(x)$ and
$$
X_t=X_0+\int^t_0 b(s, X_s,\rho_s(X_s))\dif s+\sqrt{2}W_t,\ \ \bP-a.s.
$$
\end{enumerate}
\ed

To show the existence of a weak solution, we consider the following Euler scheme to DDSDE \eqref{SDE}:
Let $T>0$, $N\in\mN$ and $h:=T/N$. For $t\in[0,h]$, define
$$
X^{N}_t:=X_0+\sqrt{2}W_{t},
$$
and for $t\in[kh,(k+1)h]$ with $k=1,\cdots,N$, we inductively define $X^N_t$ by
\begin{align}\label{SA}
X^N_t:=X^N_{kh}+(t-kh)b(kh,X^N_{kh},\rho^N_{kh}(X^N_{kh}))+\sqrt{2}(W_t-W_{kh}),
\end{align}
where $\rho^N_{kh}(x)$ is the distributional density of $X^N_{kh}$, whose existence is easily seen from the construction.
We have the following existence and uniqueness result.
\bt\label{MRE}
Assume that $b$ is bounded measurable and 
 \begin{align}\label{PR2}
\lim_{(t,u)\to(t_0,u_0)}\sup_{|x|<R}|b(t,x,u)-b(t_0,x,u_0)|=0,\ \ \forall R>0.
\end{align}
{\bf (Existence)} For any $T>0$ and initial distribution $\nu_0$,
there are a subsequence $N_k$ and a weak solution $(X,W)$ to DDSDE \eqref{SDE} in the sense of Definition \ref{Def1} 
so that for any bounded measurable $f$ and $t\in(0,T]$,
$$
\lim_{k\to\infty}\mE f(X^{N_k}_t)=\bE f(X_t).
$$
Moreover, for each $t\in(0,T]$, $X_t$ admits a density $\rho_t$ satisfying the estimate
$$
\rho_t(y)\leq Ct^{-d/2}\int_{\mR^d}\e^{-\frac{|x-y|^2}{\lambda t}}\nu_0(\dif x),\ y\in\mR^d,
$$
where $C,\lambda\geq 1$, and the following $L^1$-convergence holds:
\begin{align}\label{Li}
\lim_{k\to\infty}\int_{\mR^d}|\rho^{N_k}_t(y)-\rho_t(y)|\dif y=0.
\end{align}
{\bf (Uniqueness)} Suppose that $\nu_0(\dif x)=\rho_0(x)\dif x$ with $\rho_0\in(L^1\cap L^q)(\mR^d)$ for some $q\in(d,\infty]$, and there is a  $C>0$ such that for all $t,x,u,u'$,
\begin{align}\label{Lip}
|b(t,x,u)-b(t,x,u')|\le C|u-u'|. 
\end{align}
Then weak and strong uniqueness hold for SDE \eqref{SDE}. 
\et
\br\rm
We emphasize that the continuity of $b$ in the time variable is not necessary for the existence of weak solutions of \eqref{SDE}.
Here we need it because we are considering the Euler scheme.
Moreover, if the uniqueness holds, then limit \eqref{Li} holds for the whole sequence.
\er

As a consequence of Theorem \ref{MRE}, we have the following well-posedness of the nonlinear FPE \eqref{P13}.
\bc
Let $\nu_0$ be a probability measure over $\mR^d$. 
\begin{enumerate}[(i)]
\item Assume $b$ is bounded and Borel measurable such that \eqref{PR2} holds. Then there is a  weak solution 
$\rho_t$ to PDE \eqref{P13} in the sense  \eqref{DPDE} with $\int \rho_t=1$ and
\begin{align}\label{DA1}
0\leq \rho_t(x)\leq Ct^{-d/2}\int_{\mR^d}\e^{-\frac{|x-y|^2}{\lambda t}}\nu_0(\dif y),\ x\in\mR^d, \ t\in(0,T].
\end{align}
\item Assume that \eqref{Lip} holds and that $\nu_0(\dif x)=\rho_0(x)\dif x$ with $\rho_0\in(L^1\cap L^q)(\mR^d)$ for some $q\in(d,\infty]$.
Then the solution in assertion (i) is the unique weak solution to \eqref{P13}.
\end{enumerate}
\ec
\begin{proof}
(i) The existence of a weak solution follows from the existence part of Theorem \ref{MRE} and It\^o's formula.
\medskip\\
(ii) The uniqueness follows from the uniqueness part of Theorem \ref{MRE} by Section 2 in \cite{Ba-Ro18}. 
\end{proof}

McKean-Vlasov SDE of Nemytskii-type (density dependent SDE), i.e.
\begin{align*}
&\dif X_t=b\Big(t,X_t,\tfrac{\dif \cL_{X_t}}{\dif x}(X_t)\Big)\dif t+\sigma\Big(t,X_t,\tfrac{\dif \cL_{X_t}}{\dif x}(X_t)\Big)\dif W_t,\ \ 
X_0\stackrel{(d)}{=}\nu_0,
\end{align*}
were first introduced in \cite{Ba-Ro18a}*{Section 2}. In \cite{Ba-Ro18a}, \cite{Ba-Ro18}, 
for a large class of time independent coefficients $b,\sigma$,
 Barbu together with the second named author
obtained the existence of weak solutions for  such (possibly degenerate) density dependent SDEs (see Section 2). The strategy in   \cite{Ba-Ro18a} or \cite{Ba-Ro18}
is to solve the associated  nonlinear Fokker-Planck equation and then by the well-known superposition principle (cf. \cite{Tr}, generalizing \cite{Ku} and \cite{Fi}) 
to establish the existence of a weak solution to DDSDE \eqref{SDE}.
Later, in \cite{Ba-Ro19}, the same authors prove the uniqueness of the weak solution to DDSDE \eqref{SDE}, which is a consequence of the uniqueness of  the corresponding nonlinear Fokker-Planck equation. 
Recently, in \cite{Ba-Ro20}, they also consider the existence of solutions to a class of nonlinear Fokker-Planck equations with measure-valued initial data.
The strategy of this paper is completely different from that in \cite{Ba-Ro18a,Ba-Ro18,Ba-Ro19,Ba-Ro20}. We start directly from DDSDE \eqref{SDE} and prove weak existence of solutions 
through Euler's approximation
for DDSDE \eqref{SDE} and using simple heat kernel estimates. In other words, we do not use the superposition principle.
Moreover, our assumptions on the drift are weaker. Especially, there is no regularity assumption of $b$ in $x$. On the other hand we only consider the 
case $\sigma=$ identity.

\medskip

As explained above we obtain weak (i.e. in the sense of Schwartz distributions) solutions to FPE \eqref{P13}. Let us mention here that Chen and Perthame in \cite{CP} studied the Cauchy problem for a general nonlinear degenerate parabolic-hyperbolic equation of second order
in the framework of kinetic and entropy solutions, in the case the coefficient do not depend explicitly on $x$ and $t$. 

\medskip

However, all the above results do not cover Theorem \ref{MRE} above. 
In particular, we use a purely probabilistic method to show the existence of weak solutions for the nonlinear FPE \eqref{P13}.
This is the main contribution of the present paper. Here, for simplicity, we only consider the additive noise case. For the case of uniformly elliptic diffusion coefficients,
it would also work by the corresponding heat kernel estimates (see \cite{CHXZ} and \cite{MPZ}).

\medskip

This paper is organized as follows: In Section 2, we give some necessary preliminaries about heat kernel estimates for the Euler scheme with bounded measurable drifts. 
In Section 3, we prove our main Theorem \ref{MRE} by Euler's type approximation (cf. \cite{Z19}). 
Note that the usual Picard iteration does not seem to work for DDSDE \eqref{SDE},  since $x\to b(t,x,\rho_t(x))$ is too singular.

\medskip

Throughout this paper, we use the following conventions: The letter $C$ denotes a constant, whose value may change in different places.
We also use $A\lesssim B$  to denote $A\leq C B$ for some unimportant constant $C>0$.

\section{Heat kernel of Euler scheme}

In this section we show heat kernel estimates for Euler's scheme of usual SDEs.
First of all, we recall some basic properties about the Gaussian heat kernel. Let
\begin{align}\label{GG0}
g(t,x):=(4\pi t)^{-\frac d2}\e^{-\frac{|x|^2}{4t}},\ t>0,x\in\mR^d,
\end{align}
which is the fundamental solution of $\Delta$, i.e.,
$$
\p_t g(t,x)=\Delta g(t,x).
$$
Moreover, we have the following Chapman-Kolmogorov equations:
\begin{align}\label{CK}
(g(t)*g(s))(x):=\int_{\mR^d} g(t,x-z)g(s,z)\dif z=g(t+s,x),\ t,s>0,
\end{align}
and the following easy facts,
\begin{align}\label{GG1}
g(t,x+y)\leq 2^{\frac d2}g(2t,x)\e^{\frac{|y|^2}{4t}},\ \ |\nabla g|(t,x)\leq \tfrac{2^{d/2}}{\sqrt{t}}g(2t,x).
\end{align}

The following lemma is straightforward and elementary. For the readers' convenience, we provide a detail proof.
\bl\label{L001}
For any $T>0$, $\beta\in(0,1)$ and $j=0,1$, there is a constant  $C=C(T,\beta,j,d)>0$ such that 
for any $0<t\le T$ and $x_1,x_2\in\mR^d$,
\begin{align}\label{GS}
|\nabla^j g(t,x_1)-\nabla^j g(t,x_2)|\le C|x_1-x_2|^{\beta}t^{-\frac j2-\beta}\sum_{i=1,2}g(4 t,x_i),
\end{align}
and for any $0<t_1<t_2\le T$ and $x\in\mR^d$,
\begin{align}\label{GT}
|\nabla^j g(t_1,x)-\nabla^j g(t_2,x)|\le C|t_2-t_1|^{\frac\beta 2}\sum_{i=1,2} t_i^{-\frac{j+\beta}2}g(2 t_i,x).
\end{align}
\el
\begin{proof}
(i) By definition \eqref{GG0}, it is easy to see that for $k=1,2,3$, there is a constant $C>0$ only depending on $k,d$ such that
$$
|\nabla^k g(t,x)|\leq C (4\pi t)^{-\frac d2}t^{\frac k2}\e^{-\frac{|x|^2}{8 t}}=C2^{\frac d2} t^{-\frac k2} g(2 t,x).
$$
Thus, for $j=0,1$ and $\beta\in(0,1)$,  if $|x_1-x_2|>\sqrt{t}$, then
\begin{align*}
&|\nabla^j g(t,x_1)-\nabla^j g(t,x_2)|\lesssim t^{-\frac j2}(g(2 t,x_1)+g(2 t,x_2))
\\&\qquad\lesssim |x_1-x_2|^{\beta}t^{-\frac j2-\beta}(g(2 t,x_1)+g(2 t,x_2));
\end{align*}
if $|x_1-x_2|\leq\sqrt{t}$, then by the mean-value formula,
\begin{align*}
&|\nabla^j g(t,x_1)-\nabla^j g(t,x_2)|
\leq|x_1-x_2|\int^1_0|\nabla^{j+1} g(t,x_1+\theta(x_2-x_1))|\dif\theta
\\&\qquad\qquad\lesssim |x_1-x_2| t^{-\frac{j+1}2}\int^1_0g(2 t,x_1+\theta(x_2-x_1))\dif\theta
\\&\qquad\qquad\lesssim |x_1-x_2|t^{-\frac{j+1}2}g(4t,x_1)\lesssim |x_1-x_2|^{\beta} t^{-j/2-\beta} g(4 t,x_1).
\end{align*}
Combining the above calculations, we get \eqref{GS}.

(ii) If $t_2-t_1\leq t_1$, then by the mean-value formula,
\begin{align*}
&|\nabla^j g(t_1,x)-\nabla^j g(t_2,x)|
\leq|t_1-t_2|\int^1_0|\nabla^j\p_tg|(t_1+\theta(t_2-t_1),x)\dif\theta
\\&\qquad\qquad= |t_1-t_2|\int^1_0|\nabla^j\Delta g|(t_1+\theta(t_2-t_1),x)\dif\theta
\\&\qquad\qquad\lesssim |t_1-t_2|\int^1_0\frac{g(2 (t_1+\theta(t_2-t_1)),x)}{(t_1+\theta(t_2-t_1))^{1+j/2}}\dif\theta
\\&\qquad\qquad\lesssim |t_1-t_2|t_1^{-1-\frac j2} g(2 t_2,x)\lesssim |t_1-t_2|^{\frac \beta2}t_2^{-\frac\beta2} g(2 t_2,x);
\end{align*}
if $t_2-t_1>t_1$, then $t_2\leq 2(t_2-t_1)$ and
\begin{align*}
&|\nabla^j g(t_1,x)-\nabla^j g(t_2,x)|
\lesssim t_1^{-\frac j2}g(2 t_1,x)+t_2^{-\frac j2}g(2 t_2,x)
\\&\quad\lesssim |t_1-t_2|^{\frac\beta2}\Big(t_1^{-\frac{j+\beta}2} g(2 t_1,x)+t_2^{-\frac{j+\beta}2} g(2t_2,x)\Big).
\end{align*}
The proof is complete.
\end{proof}

Let  $b:\mR_+\times\mR^d\to\mR^d$ be a bound measurable function. Fix $T>0$ and $x\in\mR^d$. 
Let $X^N_t=X^N_t(x)$ be defined by the following Euler scheme:
\begin{align}\label{00}
X^N_t=x+\int_0^tb(\phi_N(s),X^N_{\phi_N(s)})\dif s+\sqrt{2}W_t,\quad t\in[0,T],
\end{align}
where $\phi_N(s):=kT/N$ for $s\in[kT/N,(k+1)T/N)$. 
We have the following Duhamel formula.
\bl\label{Le22}
For each $t\in(0,T]$ and $x\in\mR^d$, $X^N_t(x)$ admits a density $p^N_x(t,y)$ which satisfies the following Duhamel formula:
\begin{align}\label{Duh}
p^N_x(t,y)=g(t,x-y)+\int^t_0\mE \Big[b(\phi_N(s),X^N_{\phi_N(s)})\cdot\nabla g(t-s,y-X^N_s)\Big]\dif s.
\end{align}
\el
\begin{proof}
Fix $t\in(0,T]$ and $f\in C^\infty_c(\mR^d)$. For $s\in[0,t]$, let $u(s,x):=g(t-s,\cdot)*f(x)$. 
Since $(\p_s+\Delta)u\equiv 0$ and $u(t,x)=f(x)$, by It\^o's formula, we have 
\begin{align*}
\mE f(X^N_t)=\mE u(t, X^N_t)=u(0, x)+\int^t_0\mE\Big[b(\phi_N(s),X^N_{\phi_N(s)})\cdot\nabla u(s,X^N_s)\Big]\dif s.
\end{align*}
From this, we derive the desired Duhamel formula.
\end{proof}

\br\rm
For a general initial value $X^N_0=X_0\in\sF_0$ and each $t\in(0,T]$, since for each $x\in\mR^d$, $X^N_t(x)$ is independent of $X_0$,
the Euler scheme $X^N_t$ defined by \eqref{00} with initial value $X_0$
also has a density $p^N_{X_0}(t,y)$ given by
\begin{align}\label{Com}
p^N_{X_0}(t,y)=\int_{\mR^d}p^N_x(t,y)\mP\circ X^{-1}_0(\dif x).
\end{align}
\er

The following Gaussian type estimate for $p^N_x(t,y)$ was proved by Lemaire and Menozzi \cite{Le-Me}.
Since it is not  difficult, for the readers' convenience, we provide a detailed proof here.
\bt\label{EHE}
For any $T>0$, there is a constant $C=C(d,T,\|b\|_{\infty})$ such that for all $N\in\mN$, $t\in(0,T]$ and $x,y\in\mR^d$,
\begin{align}\label{ES09}
p^N_x(t,y)\le Cg(4t,x-y).
\end{align}
\et

\bpf
Let $\eps>0$ be small enough so that 
$$
\ell_\eps:=2^{d+1}\sqrt{\eps\|b\|_\infty^2}\e^{\eps\|b\|_\infty^2}\leq 1/2.
$$
Fix $T>0$. Without loss of generality, we assume 
\begin{align}\label{WW1}
N\geq(\|b\|^2_\infty T/(4\log 2))\vee (T/\eps).
\end{align}
For simplicity we shall write 
$$
h:=T/N,\ \ M:=[\eps/h]\in\mN.
$$
{\it Step 1:} 
In this step we use induction to show that for all $k=1,\cdots, M\wedge N$,
\begin{align}\label{Step1}
p^N_x(kh,y)\le 2^{d+1}g(4kh,x-y).
\end{align}
First of all, for $k=1$, since $X^N_h=x+W_h+hb(0,x)$, by \eqref{GG1} and \eqref{WW1} we have
\begin{align*}
p^N_x(h,y)&=g(h,y-b(0,x)h-x)\le  2^{d/2}\e^{\|b\|_\infty^2 h/4}g(2h,x-y)
\\&\leq 2^d\e^{\|b\|_\infty^2 T/(4N)}g(4h,x-y)\leq 2^{d+1}g(4h,x-y).
\end{align*}
Suppose now that \eqref{Step1} holds for $j=1,2,..,k-1$.
By Duhamel's formula \eqref{Duh}, we have
 \begin{align}\label{1003}
&p^N_x(kh,y)-g(kh,x-y)
=\int_0^{kh}\mE\[b(\phi_N(s),X^N_{\phi_N(s)})\cdot\nabla g(kh-s,y-X^N_s)\]\dif s\no\\
&\qquad\qquad=\sum_{j=0}^{k-1}\int_{jh}^{(j+1)h}\mE\[b(jh,X^N_{jh})\cdot\nabla g(kh-s,y-X^N_s)\]\dif s.
\end{align}
Note that for $s\in(jh,(j+1)h)$,
\begin{align*}
X^N_s=X^N_{jh}+\sqrt{2}(W_s-W_{jh})+(s-jh)b(jh,X^N_{jh}).
\end{align*}
Since $\sqrt{2}(W_s-W_{jh})$ is independent of $X^N_{jh}$ and has density $g(s-jh,y)$, 
by the C-K equations \eqref{CK} we have
\begin{align*}
\sI_j(s)&:=\mE\[b(jh,X^N_{jh})\cdot\nabla g(kh-s,X^N_s-y)\]
\\&=\mE\[b(jh,X^N_{jh})\cdot\nabla g(kh-s)*g(s-jh)\Big(X^N_{jh}+(s-jh)b(jh,X^N_{jh})-y\Big)\]
\\&=\mE\[b(jh,X^N_{jh})\cdot\nabla g\Big(kh-jh,X^N_{jh}+(s-jh)b(jh,X^N_{jh})-y\Big)\]
\\&\leq \|b\|_\infty\int_{\mR^d}|\nabla g|\Big(kh-jh,z+(s-jh)b(jh,z)-y\Big)p^N_x(jh,z)\dif z.
\end{align*}
By \eqref{GG1} and induction hypothesis, we further have for $s\in(jh,(j+1)h)$,
\begin{align*}
\sI_j(s)&\leq \frac{\|b\|_\infty 2^{d/2}}{\sqrt{kh-jh}}\int_{\mR^d}g\Big(2(kh-jh),z+(s-jh)b(jh,z)-y\Big)p^N_x(jh,z)\dif z
\\&\leq \frac{\|b\|_\infty 2^{d}\e^{(k-j)h\|b\|_\infty^2/4}}{\sqrt{kh-jh}}\int_{\mR^d}g(4(kh-jh),z-y)\cdot 2^{d+1} g(4jh,x-z)\dif z
\\&\leq\frac{\|b\|_\infty 2^{2d+1}\e^{\eps\|b\|_\infty^2/4}}{\sqrt{kh-s}}g(4kh,x-y)
=\frac{2^d\ell_\eps/\sqrt{\eps}}{\sqrt{kh-s}}g(4kh,x-y),
\end{align*}
where we have used $kh\leq Mh\leq\eps$.
Substituting this into \eqref{1003}, we obtain
\begin{align*}
|p^N_x(kh,y)-g(kh,x-y)|&\le \frac{2^d\ell_\eps}{\sqrt{\eps}} g(4kh,x-y)\sum_{j=0}^{k-1}\int_{jh}^{(j+1)h}\frac{1}{\sqrt{kh-s}}\dif s\\
&\le \frac{2^d\ell_\eps}{\sqrt{\eps}} g(4kh,x-y)\int^{kh}_0\frac{1}{\sqrt{kh-s}}\dif s\\
&=\frac{2^d\ell_\eps}{\sqrt{\eps}}g(4kh,x-y)2\sqrt{kh}\leq 2^{d+1}\ell_\eps g(4kh,x-y),
\end{align*}
which implies, since $g(t,x)\leq 2^d g(4t,x)$ and $2\ell_\eps\leq 1$, that
 \begin{align*}
p^N_x(kh,y)\le 2^d(1+2\ell_\eps)g(4kh,x-y)\le 2^{d+1} g(4kh,x-y).
\end{align*}
{\it Step 2:} Next we assume $M<N$ and consider $k=M+1\cdots, 2M$. Note that
\begin{align*}
X^N_{t+Mh}&=X^N_{Mh}+W_{t+Mh}-W_{Mh}+\int^{t+Mh}_{Mh}b(\phi_N(s), X^N_{\phi_N(s)})\dif s\\
&=X^N_{Mh}+W_{t+Mh}-W_{Mh}+\int^t_0b(\phi_N(s)+Mh, X^N_{\phi_N(s)+Mh})\dif s,
\end{align*}
where we have used that $\phi_N(s+Mh)=\phi_N(s)+Mh$.
In particular, if we let 
$$
\bar{X}^N_{t}:={X}^N_{t+Mh},\ \ \bar{W}_{t}:=W_{t+Mh}-W_{Mh},
$$
then for $t\in[0,Mh]$,
\begin{align*}
\bar{X}^N_{t}=X^N_{Mh}+\bar{W}_{t}+\int_0^tb({\phi}_N(s)+nMh,\bar{X}^N_{\phi_N(s)})\dif s.
\end{align*}
Let $\bar p^N_x(kh,y)$ be the density of $\bar{X}^N_{t}$ with $\bar{X}^N_0=x$.
By Step 1, we have
\begin{align*}
\bar p^N_x(kh,y)\leq 2^{d+1}g(4kh,x-y),\ k=1,\cdots, M.
\end{align*}
Thus, for $k=1,\cdots, M$, by \eqref{Com} we have
\begin{align*}
p^N_x((k+M)h,y)&=\int_{\mR^d}\bar p^N_z(kh,y)p^N_x(Mh,z)\dif z\\
&\leq 4^{d+1}\int_{\mR^d}g(4kh,z-y)g(4Mh,x-z)\dif z\\
&=4^{d+1}g(4(k+M)h,x-y).
\end{align*}
Repeating the above procedure $[\frac T\eps]+1$-times, we obtain that for some $C>0$ independent of $N$,
$$
p^N_x(kh,y)\leq Cg(kh,x-y),\ \ k=1,\cdots, N.
$$
{\it Step 3:} Note that for $t\in(kh,(k+1)h)$,
\begin{align*}
X^N_t=X^N_{kh}+W_t-W_{kh}+(t-kh)b(kh,X^N_{kh}),
\end{align*}
where $W_t-W_{kh}$ is independent of $X^N_{kh}$. Hence, 
\begin{align*}
p^N_x(t,y)&=\int_{\mR^d}g(t-kh, z+(t-kh)b(kh,z)-y)p^N_x(kh,z)\dif z\\
&\leq C\e^{(t-kh)\|b\|^2_\infty/4}\int_{\mR^d}g(4(t-kh),y-z)g(4kh,x-z)\dif z\\
&\leq C\e^{T\|b\|^2_\infty/4}g(4t,x-y).
\end{align*}
This completes the proof.
\epf

The following corollary is a combination of Theorem \ref{EHE} and Lemma \ref{L001}.
\bc\label{Cor25}
Let $\nu_0(\dif y)=\mP\circ X_0^{-1}(\dif y)$ be the distribution of $X_0$. 
\begin{itemize}
\item[(i)] For any $T>0$, there is a constant $C=C(d,T,\|b\|_\infty)$ such that for all $N\in\mN$, $t\in(0,T)$ and $x\in\mR^d$
\begin{align}\label{901}
p^N_{X_0}(t,y)\le C\int_{\mR^d} g(4t,x-y)\nu_0(\dif x).
\end{align}
\item[(ii)]
For any $T>0$ and $\beta\in(0,1)$, there is a constant $C=C(d,T,\|b\|_\infty,\beta)$ 
such that for all $N\in\mN$, $t\in(0,T)$ and $x_1,x_2\in\mR^d$,
\begin{align*}
|p^N_{X_0}(t,x_1)-p^N_{X_0}(t,x_2)|\le C|x_1-x_2|^\beta t^{-\frac\beta2}\sum_{j=1,2}\int_{\mR^d} g(4t,x_j-y)\nu_0(\dif y).
\end{align*}
\item[(iii)]
For any $T>0$ and $\beta\in(0,1)$, there is a constant $C=C(d,T,\|b\|_\infty,\beta)$ 
such that for all $N\in\mN$, $t_1,t_2\in(0,T)$ and $x\in\mR^d$,
\begin{align*}
|p^N_{X_0}(t_1,x)-p^N_{X_0}(t_2,x)|\le C|t_1-t_2|^{\beta/2}\sum_{j=1,2}t_j^{-\beta/2}\int_{\mR^d} g(2t_j,x-y)\nu_0(\dif y).
\end{align*}

\end{itemize}

\ec
\bpf
(i) is a direct consequence of \eqref{Com} and Theorem \ref{EHE}.
We only show (iii) since (ii) is similar by \eqref{GS}. Suppose $t_1<t_2$. By \eqref{Duh}, we have
\begin{align*}
&|p^N_{X_0}(t_1,y)-p^N_{X_0}(t_2,y)|\le \int_{\mR^d}| g(t_1,x-y)- g(t_2,x-y)|\nu_0(\dif x)\\
&\quad+\|b\|_{\infty}\int_0^{t_1}\!\!\int_{\mR^d}|\nabla g(t_1-s,y-z)-\nabla g(t_2-s,y-z)|p^N_{X_0}(s,z)\dif z\dif s\\
&\quad+\|b\|_\infty\int_{t_1}^{t_2}\!\!\int_{\mR^d}|\nabla g(t_1-s,y-z)|p^N_{X_0}(s,z)\dif z\dif s=:I_1+I_2+I_3.
\end{align*}
For $I_1$, by \eqref{GT}, we have
$$
I_1\lesssim |t_1-t_2|^{\frac\beta2}\sum_{j=1,2}t_j^{-\frac\beta 2} \int_{\mR^d}
g(2t_j,x-y)\nu_0(\dif y).
$$
For $I_2$, by  (i), \eqref{GT} and the C-K equations \eqref{CK}, we have
\begin{align*}
I_2&\lesssim |t_1-t_2|^{\frac\beta2}\sum_{j=1,2}\int_0^{t_1}\Bigg[\int_{\mR^d}(t_j-s)^{-\frac{1+\beta} 2} g(4(t_j-s),x-y)
\\&\qquad\qquad\qquad\qquad\qquad\times\int_{\mR^d}g(4s,z-y)\nu_0(\dif z)\dif y\Bigg]\dif s\\
&=|t_1-t_2|^{\frac\beta2}\sum_{j=1,2}\int_0^{t_1}(t_j-s)^{-\frac{1+\beta} 2} \int_{\mR^d}g(4t_j,z-x)\nu_0(\dif z)\dif s
\\&\lesssim |t_1-t_2|^{\frac\beta2}\sum_{j=1,2}\int_{\mR^d}g(4t_j,z-x)\nu_0(\dif z).
\end{align*}
For $I_3$, by  (i), \eqref{GT} and the C-K equations, we have
\begin{align*}
I_3&\lesssim\int_{t_1}^{t_2}\!\!\int_{\mR^d} (t_1-s)^{-\frac12}g(4(t_1-s),x-y)\int_{\mR^d} g(4s,y-z)\nu_0(\dif z)\dif y\dif s\\
&=\int_{t_1}^{t_2}\!\!\int_{\mR^d} (t_1-s)^{-\frac12}g(4t_1,x-z)\nu_0(\dif z)\dif s
\lesssim|t_2-t_1|^{\frac12}\int_{\mR^d}g(4t_1,x-z)\nu_0(\dif z).
\end{align*}
Combining the above calculations, we obtain the desired estimate.
\epf

\section{Proof of Theorem \ref{MRE}}

 Let $(\Omega,\sF,(\sF_t)_{t\ge0},\mP)$ be a complete filtered probability space,
 $W_t$ a $d$-dimensional standard $\sF_t$-Brownian motion, and
 $X_0$ an $\sF_0$-measurable random variable with distribution $\nu_0$.
Let $T>0$, $N\in\mN$ and $h:=T/N$. 
Let $X^N_t$ be the Euler approximation of DDSDE \eqref{SDE} constructed in the introduction.
From the construction \eqref{SA}, it is easy to see that $X^N_t$ solves the following SDE:
\begin{align}\label{SDEN}
X_t^N=X_0+\int_0^tb^N(\phi_N(s),X_{\phi_N(s)}^N)\dif s+\sqrt{2}W_t,
\end{align}
where
\begin{align}\label{sb}
b^N(s,x)=\1_{\{s\ge h\}}b\Big(s,x,\rho^N_{\phi_N(s)}(x)\Big)
\end{align}
and
\begin{align}\label{sb0}
\phi_N(s):=\sum_{j=0}^\infty jh\1_{[jh,(j+1)h)}(s).
\end{align}
The following lemma is easy by \eqref{SDEN} and since $\|b^N\|_\infty\leq\|b\|_\infty$.
\bl\label{L03}
For any $T>0$, there is a constant $C>0$ such that for all
$s,t\in[0,T]$,
$$
\sup_{N}\mE|X_t^N-X_s^N|^4\le C|s-t|^2.
$$
\el

Let $p^N_x(t,y)$ be the distributional density of the Euler scheme $X^N_t(x)$  of SDE \eqref{SDEN} starting from $x$ at time $0$. 
Since for each $x\in\mR^d$, $X^N_t(x)$ is independent of $X_0$, the distributional density $\rho^N_t(y)$ of $X^N_t$ with 
initial distribution $\nu_0$ is given by
\begin{align}\label{PR0}
\rho^N_t(y)=\int_{\mR^d}p^N_x(t,y)\nu_0(\dif x).
\end{align}

The following lemma is crucial for the existence of a solution to DDSDE \eqref{SDE}.
\bl\label{AA}
For fixed $T>0$, there are a subsequence $(N_k)_{k\in\mN}$ and a continuous function $\rho\in C((0,T]\times\mR^d)$ such that for any $M\in\mN$,
\begin{align}\label{USC}
\lim_{k\to\infty}\sup_{|y|\le M}\sup_{1/M\le t\le T}|\rho^{N_k}_t(y)-\rho_t(y)|=0.
\end{align}
\el

\bpf
First of all, by the upper-bound estimate \eqref{ES09} for $p^N_x(t,y)$, we have
$$
 \sup_{|y|\le M}\sup_{1/M\le t\le T}|\rho^{N}_t(y)|\le C\int_{\mR^d}\sup_{|y|\le M}\sup_{1/M\le t\le T}| g(4t,x-y)|\nu_0(\dif x)\le C_M,
$$
where $C_M$ is independent of $N$.
Moreover, by Corollary \ref{Cor25}, we have for any $\beta<1$, $t_1,t_2\in[1/M,T]$ and $y_1,y_2\in\mR^d$,
\begin{align}
|\rho^{N}_{t_1}(y_1)-\rho^{N}_{t_2}(y_2)|&\le|\rho^{N}_{t_1}(y_1)-\rho^{N}_{t_2}(y_1)|+|\rho^{N}_{t_2}(y_1)-\rho^{N}_{t_2}(y_2)|\no\\
&\lesssim|t_2-t_1|^{\frac\beta 2}\sum_{i=1,2}\int_{\mR^d}| g(2t_i,y_1-x)|\nu_0(\dif x)\no\\
&+|y_1-y_2|^{\beta}\sum_{i=1,2}\int_{\mR^d}| g(4t_2,y_i-x)|\nu_0(\dif x)\no\\
&\lesssim  M^{-(d+1+\beta)/2}\(|t_2-t_1|^{\frac\beta 2}+|y_1-y_2|^{\beta}\),\label{PR1}
\end{align}
where the implicit constants in the above $\lesssim$ are independent of $N$.
Thus, by Ascolli-Arzela's theorem, we conclude the proof and have \eqref{USC}.
\epf

Now we are in a position to give {\red} the

\begin{proof}[Proof of \autoref{MRE}]
{\bf (Existence)} Fix $T>0$. Let $\mW$ be the space of all continuous functions from $[0,T]$ to $\mR^d$.
 Let $\mQ_N$ be the law of $(X_\cdot^N,W_\cdot)$ in $\mW\times\mW$. By Lemma \ref{L03} and Kolmogorov's criterion, 
 $\{\mQ_N\}_{N\in\mN}$ is tight. Therefore, by Prokhorov's theorem, there are 
 a subsequence $(N_k)_{k\in\mN}$ and a probability measure $\mQ$ on $\mW\times\mW$ so that
$$
\mQ_{N_k}\to\mQ\quad\text{weakly.}
$$
Without loss of generality, we assume that the subsequence is the same as that in Lemma \ref{AA}.
Below, we still denote the above subsequence by $\mQ_N,~N\in\mN$ for simplicity.
Now, by Skorokhod's representation theorem, there are probability space $(\widetilde {\Omega},\widetilde {\sF},\widetilde {\mP})$ 
and random variables $(\widetilde {X}^N,\widetilde {W}^N)$ and $(\widetilde {X},\widetilde {W})$ thereon such that
\begin{align}\label{Ea.e}
(\widetilde {X}^N,\widetilde {W}^N)\to(\widetilde {X},\widetilde {W}),\quad\text{$\widetilde {\mP}$-a.s.}
\end{align}
and
\begin{align}\label{ID}
\widetilde {\mP}\circ(\widetilde {X}^N,\widetilde {W}^N)^{-1}=\mQ_N={\mP}\circ({X}^N,{W})^{-1},\quad \widetilde {\mP}\circ(\widetilde {X},\widetilde {W})^{-1}=\mQ.
\end{align}
In particular, the distributional density of $\widetilde {X_t^N}$ is $\rho^N_t$. Moreover, 
by Lemma \ref{AA} and \eqref{Ea.e}, for any $t\in(0,T)$ and $\varphi\in C_0^\infty(\mR^d)$,
\begin{align*}
\mE\varphi(\widetilde {X}_t)=\lim_{N\to\infty}\mE\varphi(\widetilde {X}^N_t)=\lim_{N\to\infty}\int_{\mR^d}\varphi(y)\rho^N_t(y)\dif y=\int_{\mR^d}\varphi(y)\rho_t(y)\dif y.
\end{align*} 
In other words, $\rho_t$ is the density of $\widetilde {X}_t$.
Define $\widetilde \sF_t^N:=\sigma(\widetilde {X}^N,\widetilde {W}^N;s\le t)$. We note that
$$
\mP[W_t-W_s\in\cdot|\sF_s]=\mP\{W_t-W_s\in\cdot\},
$$
hence,
$$
 \widetilde {\mP}[\widetilde W^N_t-\widetilde W^N_s\in\cdot|\widetilde \sF^N_s]=\widetilde {\mP}\{\widetilde W^N_t-\widetilde W^N_s\in\cdot\},
$$
which means that $\widetilde W_t^N$ is an $\widetilde \sF_t^N$-BM.
Thus, by \eqref{SDEN} and \eqref{ID} we have
\begin{align}\label{SDEN1}
\widetilde {X}_t^N=\widetilde {X}^N_0+\int_0^tb^N(\phi_N(s),\widetilde {X}_{\phi_N(s)}^N)\dif s+\sqrt{2}\widetilde {W}_t^N.
\end{align}
Let us now show that
\begin{align}\label{KPC}
\int_0^t \1_{s\geq h}b\Big(\phi_N(s),\widetilde {X}_{\phi_N(s)}^N,\rho^N_{\phi_N(s)}(\widetilde  X^N_{\phi_N(s)})\Big)\dif s\to\int_0^tb\Big(s,\widetilde {X}_s,\rho_s(\widetilde {X}_s)\Big)\dif s,
\end{align}
in probability as $N\to\infty$. 
Let $\Omega_0\subset\widetilde\Omega$ be a measurable set so that $\widetilde \mP(\Omega_0)=1$ and for each $\omega\in\Omega_0$,
\begin{align}\label{PR3}
\lim_{N\to\infty}(\widetilde X^N_{\phi_N(\cdot)}(\omega),\widetilde W^N_\cdot(\omega))=(\widetilde X_\cdot(\omega), \widetilde W_\cdot(\omega)).
\end{align}
In particular, for each fixed $\omega\in\Omega_0$ and $s\in(0,T)$, by \eqref{USC}, \eqref{PR1} and \eqref{PR3}, we have
\begin{align}\label{PR4}
\lim_{N\to\infty}|\rho^N_{\phi_N(s)}(\widetilde  X^N_{\phi_N(s)}(\omega))-\rho_{s}(\widetilde  X^N_{\phi_N(s)}(\omega))|=0,\ \ a.s.
\end{align}
On the other hand, for any $s>0$, by \eqref{PR2} we have
$$
\lim_{|\delta_1|+|\delta_2|\to 0}\sup_{|u|\leq R}\sup_{|x|<R}|b(s+\delta_1,x,u+\delta_2)-b(s,x,u)|=0,\ \ \forall R>0.
$$
Thus, by \eqref{PR3} and \eqref{PR4}, we have for each $s>0$ and $\omega\in\Omega_0$,
$$
\lim_{N\to\infty}|b({\phi_N(s)},\widetilde {X}_{\phi_N(s)}^N(\omega),\rho^N_{\phi_N(s)}
(\widetilde  X^N_{\phi_N(s)}(\omega)))-b(s,\widetilde {X}_{\phi_N(s)}^N(\omega),\rho_s(\widetilde  X^N_{\phi_N(s)}(\omega)))|=0,
$$
which, by  the dominated convergence theorem, implies that
\begin{align*}
&\lim_{N\to\infty}\widetilde \mE\int_0^t |\1_{s>h}b(\phi_N(s)),\widetilde {X}_{\phi_N(s)}^N,\rho^N_{\phi_N(s)}(\widetilde  X^N_{\phi_N(s)}))
-b(s,\widetilde {X}_{\phi_N(s)}^N,\rho_s(\widetilde  X^N_{\phi_N(s)}))|\dif s\\
&=\widetilde \mE\int_0^t \lim_{N\to\infty}|\1_{s>h}b(s,\widetilde {X}_{\phi_N(s)}^N,\rho^N_{\phi_N(s)}(\widetilde  X^N_{\phi_N(s)}))
-b(s,\widetilde {X}_{\phi_N(s)}^N,\rho_s(\widetilde  X^N_{\phi_N(s)}))|\dif s=0.
\end{align*}
For proving \eqref{KPC}, it remains to show
$$
\lim_{N\to\infty}\widetilde \mE\int_{h}^t |b(s,\widetilde {X}_{\phi_N(s)}^N,\rho_s(\widetilde  X^N_{\phi_N(s)}))
-b(s,\widetilde {X}_s,\rho_s(\widetilde  X_s))|\dif s=0.
$$
Let $K_\eps$ be a family of mollifiers in $\mR^d$. Define
$$
B_\eps(t,x)=b(t,\cdot,\rho_t(\cdot))*K_\eps(x).
$$
Clearly, for fixed $\eps>0$, by \eqref{PR3} we have
$$
\lim_{N\to\infty}\widetilde \mE\int_{h}^t |B_\eps(s,\widetilde {X}_{\phi_N(s)}^N)-B_\eps(s,\widetilde {X}_s)|\dif s=0.
$$
Below for notational convenience, we write $\widetilde X^\infty_t:=\widetilde X_t$ and $\phi_\infty(s):=s$.
For $N\in\mN\cup\{\infty\}$, we have
\begin{align*}
&\widetilde \mE\int_{h}^t |B_\eps(s,\widetilde {X}_{\phi_N(s)}^N)-b(s,\widetilde {X}_{\phi_N(s)}^N,\rho_s(\widetilde  X^N_{\phi_N(s)})|\dif s\\
&\leq\widetilde \mE\int_{h}^t \1_{|\widetilde  X^N_{\phi_N(s)}|\leq R}\Big|B_\eps(s,\widetilde {X}_{\phi_N(s)}^N)
-b(s,\widetilde {X}_{\phi_N(s)}^N,\rho_s(\widetilde  X^N_{\phi_N(s)})\Big|\dif s\\
&\quad+2\|b\|_\infty\int_{h}^t \widetilde\mP\Big(|\widetilde  X^N_{\phi_N(s)}|>R\Big)\dif s=:I^N_R(\eps)+J^N_R.
\end{align*}
For $I^N_R(\eps)$, by \eqref{PR0}, \eqref{ES09} and H\"older's inequality with $p>2d$ and $q=\frac{p}{p-1}$, we have
\begin{align*}
I^N_R(\eps)&=\int_{h}^t\!\int_{B_R}|B_\eps(s,y)-b(s,y,\rho_s(y))|\rho^N_{\phi_N(s)}(y)\dif y\dif s
\\&\lesssim\int_{h}^t\!\int_{B_R}|B_\eps(s,y)-b(s,y,\rho_s(y))|\int_{\mR^d}g(4\phi_N(s),x-y)\nu_0(\dif x)\dif y\dif s
\\&\lesssim\int_{h}^t\left(\int_{B_R}|B_\eps(s,y)-b(s,y,\rho_s(y))|^p\dif y\right)^{\frac 1p}
\\&\qquad\qquad\times\left(\int_{B_R}\left|\int_{\mR^d}g(4\phi_N(s),x-y)\nu_0(\dif x)\right|^q\dif y\right)^{\frac 1q}\dif s
\\&\lesssim\int_{h}^t\left(\int_{B_R}|B_\eps(s,y)-b(s,y,\rho_s(y))|^p\dif y\right)^{\frac 1p} \phi_N(s)^{-\frac dp}\dif s
\\&\lesssim\left(\int_{h}^t\left(\int_{B_R}|B_\eps(s,y)-b(s,y,\rho_s(y))|^p\dif y\right)^{\frac 2p}\dif s\right)^{\frac12}
\left(\int^t_{h} \phi_N(s)^{-\frac {2d}p}\dif s\right)^{\frac12}
\\&\lesssim\left(\int_0^t\left(\int_{B_R}|B_\eps(s,y)-b(s,y,\rho_s(y))|^p\dif y\right)^{\frac 2p}\dif s\right)^{\frac12}\left(\int^t_0 s^{-\frac {2d}p}\dif s\right)^{\frac12},
\end{align*}
where the implicit constant in the above $\lesssim$ is independent of $N, R$ and $\eps$. Hence, for each $R>0$, by the dominated convergence theorem, we obtain
$$
\lim_{\eps\to 0}\sup_{N\in\mN\cup\{\infty\}}I^N_R(\eps)=0.
$$
For $J^N_R$, by Chebyshev's inequality and \eqref{SDEN} and since $\|b^N\|_\infty\leq\|b\|_\infty$, we have
\begin{align*}
J^N_R&=2\|b\|_\infty\int_h^t\mP(|X^N_{\phi_N(s)}|>R)\dif s\\
&\leq 2\|b\|_\infty\int_0^t \mP(| X_0|+s\|b\|_\infty+\sqrt{2}|W_{\phi_N(s)}|>R)\dif s\\
&\le2\|b\|_\infty\(\int_0^t \mP(| X_0|+s\|b\|_\infty>R/2)\dif s+\int_0^t \frac{2\phi_N(s)}{(R/2)^2}\dif s\),
\end{align*}
which converges to zero uniformly in $N$, as $R\to\infty$. Combining the above calculations, we obtain
$$
\lim_{\eps\to 0}\sup_{N\in\mN\cup\{\infty\}}\widetilde \mE\int_h^t |B_\eps(s,\widetilde {X}_{\phi_N(s)}^N)
-b(s,\widetilde {X}_{\phi_N(s)}^N,\rho_s(\widetilde  X^N_{\phi_N(s)})|\dif s=0.
$$
Thus, \eqref{KPC} is proven and the existence of a solution to DDSDE \eqref{SDE} is obtained.
\medskip\\
{\bf (Uniqueness)}
Let $X_t$ and $\bar X_t$ be two solutions of DDSDE \eqref{SDE} defined on the same probability space and with the same initial value $X_0$, 
where $X_0$ has the distributional density $\rho_0\in L^q(\mR^d)$ with $q\in(d,\infty]$.
Let $\rho_t(y)$ and $\bar\rho_t(y)$ be the distributional density of $X_t$ and $\bar X_t$, respectively. 
Clearly,  these are two solutions of the nonlinear Fokker-Planck equation \eqref{P13} with the same initial value $\rho_0$. 
Consider the following linearized SDE:
$$
\dif X_t=B(t,X_t)\dif t+\sqrt{2}\dif W_t,\  X_0=x,
$$
where $B(t,x):=b(t,x,\rho_t(x))$. It is well known that $X_t(x)$ admits a density $p_x(t,y)$ with Gaussian type estimate:
For some $\lambda, C>0$, it holds that for all $t\in(0,T]$ and $x,y\in\mR^d$,
$$
p_x(t,y)\leq C g(\lambda t,x-y).
$$
Note that by \eqref{ES09} and H\"older's inequality,
\begin{align}
\rho_t(y)&=\int_{\mR^d}p_x(t,y)\rho_0(x)\dif x
\lesssim\int_{\mR^d}g(\lambda t,x-y)\rho_0(x)\dif x\no\\
&\leq\|g(\lambda t,\cdot)\|_{q/(q-1)}\|\rho_0\|_q\lesssim t^{-d/(2q)}\|\rho_0\|_q.\label{DQ1}
\end{align}
Let
$$
\Gamma_t:=\rho_t-\bar\rho_t,\ \ B_t:=b(t,\cdot,\rho_t)\rho_t-b(t,\cdot,\bar\rho_t)\bar\rho_t
$$
and
$$
\Gamma^\eps_t(x):=\Gamma_t*K_\eps(x),
$$
where $\{K_\eps(x),\eps\in(0,1)\}$ is a family of mollifiers.
By definition, it is easy to see that
$$
\Gamma^\eps_t=\int^t_0\Delta\Gamma^\eps_s\dif s+\int^t_0 \div (B_s*K_\eps)\dif s.
$$
Let $\beta_\delta(r):=\sqrt{r^2+\delta}-\sqrt{\delta}$. For simplicity, we write $\int f$ for $\int_{\mR^d}f(x)\dif x$.
By the chain rule and integration by parts, we have
\begin{align*}
\p_t\int\beta_\delta(\Gamma^\eps_t)&=\int\beta_\delta'(\Gamma^\eps_t)\big(\Delta\Gamma^\eps_t+\div (B_t*K_\eps)\big)\\
&=\int\beta''_\delta(\Gamma^\eps_t)\Big(-|\nabla\Gamma^\eps_t|^2-(B_t*K_\eps)\cdot\nabla\Gamma^\eps_t\Big)\\
&\leq\int\beta''_\delta(\Gamma^\eps_t)\Big(-\tfrac{1}{2}|\nabla\Gamma^\eps_t|^2+2|B_t*K_\eps|^2\Big),
\end{align*}
where we have used that $\beta_\delta''(r)=\delta/(r^2+\delta)^{3/2}>0$.
Hence, by Fatou's lemma,
\begin{align}\label{DK9}
\int|\Gamma_t|\leq\lim_{\delta\to 0}\lim_{\eps\to0}\int\beta_\delta(\Gamma^\eps_t)
\leq 2\lim_{\delta\to 0}\lim_{\eps\to0}\int^t_0\int \beta''_\delta(\Gamma^\eps_s) (|B_s|*K_\eps)^2.
\end{align}
For fixed $\delta>0$, since $\beta''_\delta(r)\leq \delta^{-1/2}$ and by \eqref{DQ1},
$$
|B_t(x)|\leq  \|b\|_\infty(\rho_t(x)+\bar\rho_t(x))\leq C(t^{-d/(2q)}+1),
$$
the dominated convergence theorem implies
\begin{align*}
\lim_{\eps\to0}\int^t_0\!\!\int \beta''_\delta(\Gamma^\eps_s) \big|(|B_s|*K_\eps)^2-|B_s|^2\big|
&\leq \delta^{-\frac12}\int^t_0\lim_{\eps\to0}\int\big|(|B_s|*K_\eps)^2-|B_s|^2\big|=0
\end{align*}
and also,
\begin{align*}
\lim_{\eps\to0}\int^t_0\!\!\int |\beta''_\delta(\Gamma^\eps_s)-\beta''_\delta(\Gamma_s)| \cdot|B_s|^2
&=\int^t_0\!\!\int\lim_{\eps\to0} |\beta''_\delta(\Gamma^\eps_s)-\beta''_\delta(\Gamma_s)| \cdot|B_s|^2=0.
\end{align*}
Therefore, 
$$
\lim_{\eps\to0}\int^t_0\int \beta''_\delta(\Gamma^\eps_s) (|B_s|*K_\eps)^2=\int^t_0\int \beta''_\delta(\Gamma_s) |B_s|^2,
$$
and  by \eqref{DK9}, thanks to $\beta''_\delta(r)\leq1/r$,
$$
\int|\Gamma_t|\leq 2\lim_{\delta\to 0}\int^t_0\int\beta''_\delta(\Gamma_s) |B_s|^2\leq2\int^t_0\int\frac{|B_s|^2}{|\Gamma_s|}.
$$
Moreover, noting that by \eqref{Lip} and \eqref{DQ1}, 
$$
|B_t(x)|\leq  C|\Gamma_t(x)|(\|\rho_t\|_\infty+1)\leq C|\Gamma_t(x)|(t^{-\frac{d}{2q}}+1),
$$
we further have for $\gamma\in(1,\frac{q}{d})$,
\begin{align*}
\int|\Gamma_t|&\lesssim\int^t_0(s^{-\frac dq}+1)\int |\Gamma_s|
\leq\left(\int^t_0(s^{-\frac dq}+1)^\gamma \dif s\right)^{\frac1\gamma}
\left(\int^t_0\left[\int |\Gamma_s|\right]^{\frac{\gamma}{\gamma-1}}\dif s\right)^{\frac{\gamma-1}{\gamma}}.
\end{align*}
Thus by Gronwall's inequality, we get
$$
\int|\Gamma_t|\equiv0,
$$
which implies $\rho_t=\bar \rho_t.$ Now the pathwise uniqueness of SDE \eqref{SDE} follows by the well-known pathwise uniqueness for SDE \eqref{SDE} with bounded measurable drift $b(t,x,\rho_t(x))$ (cf. \cite{Ve}).
\end{proof}


\begin{bibdiv}
\begin{biblist}
\bib{Ba-Ro18a}{article}{
      author={Barbu, V.},
      author={R\"ockner, M.},
   title={Probabilistic representation for solutions to nonlinear Fokker-Planck equations},
   journal={SIAM J.Math.Anal.},
   volume={50},
	date={2018},
	number={4},
	pages={4246-4260},
   }

\bib{Ba-Ro18}{article}{
      author={Barbu, V.},
      author={R\"ockner, M.},
   title={From nonlinear Fokker-Planck equations to solutions of distribution dependent SDE},
   journal={Annals of Probability},
  volume={48},
	date={2020},
	number={4},
	pages={1902-1920},
   }
\bib{Ba-Ro19}{article}{
     author={Barbu, V.},
      author={R\"ockner, M.},
   title={Uniqueness for nonlinear Fokker-Planck equations and weak uniqueness for McKean-Vlasov SDEs},
   eprint={1909.04464},
   }


\bib{Ba-Ro20}{article}{
     author={Barbu, V.},
      author={R\"ockner, M.},
   title={Solutions for nonlinear Fokker-Planck equations with measures as initial data and McKean-Vlasov equations},
   eprint={2005.02311},
   }




\bib{CHXZ}{article}{
	author={Chen, Z.-Q.},
	author={Hu, E.},
	author={Xie, L.},
	author={Zhang, X.},
	title={Heat kernels for non-symmetric diffusion operators with jumps},
	journal={J. Differential Equations},
	volume={263},
	date={2017},
	number={10},
	pages={6576--6634},
	issn={0022-0396},
	review={\MR{3693184}},
	doi={10.1016/j.jde.2017.07.023},
}

\bib{CP}{article}{
   author={Chen, G.Q.},
   author={Perthame, B.},
   title={Well-posedness for nonisotropic degenerate parabolic-hyperbolic equations},
   journal={Ann. Inst. H. Poincar\'e (C) Non Linear Analysis},
   volume={20}
   date={2003},
   pages={645-668},
}

\bib{Fi}{article}{
   author={Figalli, A.},
   title={Existence and uniqueness of martingale solutions for SDEswith rough or degenerate coefficients},
   journal={J. Funct. Anal.},
   volume={254}
   date={2008},
   pages={109-153},
}

\bib{Le-Me}{article}{
author = {V. Lemaire}
author = {S. Menozzi},
     TITLE = {On some Non Asymptotic Bounds for the Euler Scheme},
   JOURNAL = {Electron. J. Probab.},
  FJOURNAL = {Electronic Journal of Probability},
    VOLUME = {15},
      YEAR = {2010},
     PAGES = {No. 53, 1645-1681},
}


\bib{Ku}{article}{
   author={Kurtz, T.G.},
   title={Martingale problems for conditional distributions of Markov processes},
   journal={Elec. J. Probab.},
   volume={3}
   date={1998},
   pages={1-29},
}

\bib{MPZ}{article}{
   author={Menozzi, S.},
      author={Pesce, A.},
         author={Zhang, X.},
   title={Density and gradient estimates for non degenerate Brownian SDEs with unbounded measurable drift},
      eprint={Preprint},
}
\bib{Tr}{article}{
   author={Trevisan, D.},
   title={Well-posedness of multidimensional diffusion process with weakly differentiable coefficients},
   journal={Electron.J.Probab.},
   volume={21}
   date={2016},
   number={22},
}

\bib{Ve}{article}{
   author={Veretennikov, A},
   title={On the strong solutions of stochastic differential equations},
   journal={Theory Probab. Appl.},
   volume={24}
   date={1979},
   pages={354-366},
   review={\MR{3325091}},
   doi={10.1214/EJP.v20-3287},
}

\bib{Z19}{article}{
   author={Zhang, X.},
   title={A discretized version of Krylov’s estimate and its applications},
   journal={Electron. J. Probab.},
   volume={24},
   date={2019},
   pages={1-17},
}


\end{biblist}
\end{bibdiv}
\end{document}